\documentclass{amsart}
\usepackage{amsthm}                 
\usepackage{amsmath}                
\usepackage{amsfonts,amssymb}       
\usepackage{graphicx}               
\usepackage{float}                  
\usepackage{rotating}               
\restylefloat{figure}               
\usepackage{mathrsfs}               
\usepackage{fancyhdr}               
\usepackage{hyperref}               
\usepackage[numbers]{natbib}
\usepackage{color}
\theoremstyle{plain}
\newtheorem{lem}{Lemma}

\newtheorem{thm}[lem]{Theorem}

\theoremstyle{definition}
\newtheorem{defn}[lem]{Definition}

\newtheorem{rem}[lem]{Remark}
\newcommand{\R}{\mathbb R}

\newcommand{\Z}{\mathbb Z}
\newcommand{\N}{\mathbb N}

\newcommand{\dx}{\,\text{\rm d}x}
\renewcommand{\d}{\,\text{\rm d}}
\newcommand{\dw}{\text{\rm d}}

\renewcommand{\S}{\mathbb S}
\renewcommand{\phi}{\varphi}

\newcommand{\eps}{\varepsilon}

\newcommand{\set}[2]{\left\{#1;\;#2\right\}}

\newcommand{\bea}{\begin{eqnarray}}
\newcommand{\eea}{\end{eqnarray}}
\newcommand{\beq}{\begin{equation}}
\newcommand{\eeq}{\end{equation}}
\renewcommand{\phi}{\varphi}

\renewcommand{\autoref}[1]{\text{Eq.}~\eqref{#1}}

\begin{document}
\title{On initial boundary value problems for variants of the Hunter-Saxton equation}
\author{Martin Kohlmann}
\address{Institute for Applied Mathematics, University of Hannover, D-30167 Hannover, Germany}
\email{kohlmann@ifam.uni-hannover.de}
\keywords{Hunter-Saxton equation, initial boundary value problems, well-posedness, blow-up}
\subjclass[2010]{35G31, 35B44, 35B65}
\begin{abstract}
The Hunter-Saxton equation serves as a mathematical model for orientation waves in a nematic liquid crystal.
The present paper discusses a modified variant of this equation, coming up in the study of critical points for the speed of orientation waves, as well as a two-component extension. We establish well-posedness and blow-up results for some initial boundary value problems for the modified Hunter-Saxton equation and the two-component Hunter-Saxton system.
\end{abstract}
\maketitle
\tableofcontents
\section{Introduction}
Liquid crystals are a state of matter which is intermediate between the crystalline solid and the amorphous liquid, \cite{Ch92}. For instance, a liquid crystal may flow like a liquid, but its molecules may be oriented in a crystal-like way. One of the most common phases in which a liquid crystal exists is the nematic phase. The word nematic comes from the Greek $\nu\eta\mu\alpha$ (nema), which means "thread"; in this state of matter, the molecules are oriented in a thread-like way. A simple mathematical model which neglects the kinetic energy of the molecules assums that only the orientation of the fluid particles determines the energy of the liquid crystal, cf.~\cite{HS91}. Precisely, the potential energy is obtained by restricting the Oseen-Frank potential energy function $W$ for liquid crystals to a field of unit vectors (called \emph{directors}) that lie on a circle and depend on a single space variable $x$ and time $t$:
$$\mathbf n=\cos\psi(t,x)\mathbf e_x+\sin\psi(t,x)\mathbf e_y.$$
The potential energy is now given by
\bea W(\mathbf n,\nabla\mathbf n)&=&\frac{1}{2}\left(k_1(\nabla\cdot\mathbf n)^2+k_2(\mathbf n\cdot(\nabla\times\mathbf n))^2+k_3|\mathbf n\times(\nabla\times\mathbf n)|^2\right)
\nonumber\\
&=&\frac{1}{2}\left(k_1\sin^2\psi+k_3\cos^2\psi\right)\psi_x^2,\nonumber
\eea
where the coefficients $k_1$, $k_2$ and $k_3$ are the physical parameters modeling splay, twist and bend respectively. The action functional for a nematic liquid crystal is given by
$$\mathcal S=\iint\left\{\frac{1}{2}|\mathbf n_t|^2-W(\mathbf n,\nabla\mathbf n)+\frac{1}{2}\lambda|\mathbf n|^2\right\}\dw\mathbf x\d t,$$
where $\lambda$ is a Lagrangian multiplier enforcing that $|\mathbf n|=1$; see formula (A.14) in \cite{HS91}. It follows that
$$\mathcal S=\iint\frac{1}{2}\left\{\psi_t^2-c^2(\psi)\psi_x^2\right\}\dw x\d t,\quad c^2(\psi)=k_1\sin^2\psi+k_3\cos^2\psi.$$
The corresponding Euler-Lagrange equation derived in \cite{HS91} reads
$$\psi_{tt}-c(\psi)[c(\psi)\psi_x]_x=0.$$
In search for weakly nonlinear asymptotic solutions of this equation, Hunter and Saxton suggested an expansion of the form
$$\psi(t,x)=\psi_0+\eps\psi_1(\tau,\theta)+\mathcal O(\eps^2)$$
where $\theta=x-c_0t$, $\tau=\eps t$ and $c_0=c(\psi_0)>0$ is the unperturbed wave speed. This way, one obtains the equation
$$(\psi_{1\tau}+c_0'\psi_1\psi_{1\theta})_\theta=\frac{1}{2}c_0'\psi_{1\theta}^2,\quad c_0'=c'(\psi_0),$$
which reduces, after an appropriate change of variables, to the well-known Hunter-Saxton (HS) equation
\begin{equation} (u_t+uu_x)_{x}=\frac{1}{2}u_{x}^2.\label{HS}\end{equation}
At a critical point of the wave speed (i.e., $c^{(p)}(\psi_0)\neq 0$ for some positive integer $p$ and $c^{(k)}(\psi_0)=0$ for all $0<k<p$), a suitable weakly nonlinear asymptotic expansion is of the form
$$\psi(t,x)=\psi_0+\eps\psi_1(\tau,\theta)+\mathcal O(\eps^{p+1}),\quad\theta=x-c_0t,\;\tau=\eps^pt,$$
and the resulting equation reads
$$(\psi_{1\tau}+\Gamma\psi_1^p\psi_{1\theta})_\theta=\frac{p}{2}\Gamma\psi_1^{p-1}\psi_{1\theta}^2,\quad \Gamma=\frac{c^{(p)}(\psi_0)}{p!}.$$
Writing $u=|\Gamma|^{1/p}\psi_1$, $x=(\text{sign }\Gamma)\theta$ and $t=\tau$ we find
\beq (u_t+u^pu_x)_x=\frac{p}{2}u^{p-1}u_x^2.\label{mHS}\eeq
In what follows, we will call \eqref{mHS} the modified Hunter-Saxton equation. Observe that \autoref{mHS} reduces to the Hunter-Saxton equation \eqref{HS} if $p=1$. While the Hunter-Saxton equation has been subject of a wide range of papers, its modified variant \eqref{mHS} appears only in the seminal paper \cite{HS91} and in Ti\v{g}lay's work \cite{T05} where the author establishes local well-posedness in the Sobolev spaces $H^s(\S)$ for $s>3/2$ and $p\geq 1$ and in $C^1(\S)$ for $p\geq 1$; here $\S=\R/\Z$ is the unit circle.

The mathematical theory behind the Hunter-Saxton equation is rich and interesting and represents an area of active mathematical research. Let us briefly survey some recent results: The HS equation can be regarded as a non-local perturbation of the Burgers equation. Also, it is related to the famous Camassa-Holm (CH) equation
$$m_t=-m_xu-2u_xm,\quad m=u-\alpha u_{xx},$$
from which it is formally obtained in the limit case $\alpha\to\infty$ or in the short-wave limit $(t,x)\mapsto(\eps t,\eps x)$ for $\eps\to 0$. Furthermore, the Hunter-Saxton equation belongs to the family of $b$-equations \cite{DGH04,EY08}
$$m_t=-m_xu-bu_xm,\quad m=Au,$$
from which it emerges for the choice $b=2$ and the linear operator $A=-\partial_{x}^2$. Next, the HS equation is known to re-express a geodesic flow on the diffeomorphism group of the circle (modulo rigid rotations), equipped with the right-invariant homogenous Sobolev metric induced by the $\dot H^1$ inner product on the tangent space at the identity element; cf.~\cite{L07',L07'',L08}. Its bi-Hamiltonian nature and a Lax pair representation are presented in \cite{HZ94}. An application of the inverse scattering approach to the Hunter-Saxton equation can be found in \cite{BSS01}. Another commonality of the Camassa-Holm and the Hunter-Saxton equation (and also the KdV equation) is that these equations enjoy the same symmetry group which is the Virasoro group, as explained in \cite{KM03}. Global solutions of the Hunter-Saxton equation are presented in \cite{BC05} and explicit solution formulas are calculated in \cite{L08}.

A two-component variant of the Hunter-Saxton equation (denoted as 2HS) is the subject of \cite{K10,LL09,WW10,W09,W10},
\beq
\left\{
\begin{array}{rcl}
  m_t+um_x       & = & -au_xm-\kappa\rho\rho_x, \\
  m              & = & -u_{xx},\\
  \rho_t+u\rho_x & = & (1-a)u_x\rho, \\
\end{array}
\right.
\qquad x\in\S,\;t>0,\;(a,\kappa)\in\R^2,
\label{2HS}\eeq
where the authors discuss geometric aspects, integrability and well-posedness issues. Observe that for $\rho\equiv 0$ (which can be achieved by requiring $\rho_0\equiv 0$) and $a=2$, the two-component Hunter-Saxton equation reduces to the one-component version \autoref{HS}. Some of the results in \cite{WW10,W10} about the two-component Hunter-Saxton system \eqref{2HS} will be important for this paper.\\

\textbf{Outline of results.} The first aim of the paper at hand is to study two kinds of initial boundary value problems for the modified Hunter-Saxton equation on the interval $(0,\frac{1}{2})$, with initial data in $H^s(0,\frac{1}{2})$ and $D^s_k(0,\frac{1}{2})$ respectively, where the latter space is the subspace of $H^s(0,\frac{1}{2})$ which consists of functions $v$ for which
$$v^{(2k-2j)}(0)=v^{(2k-2j)}(\tfrac{1}{2})=0,\quad j=0,\ldots,k,\,k\in\N,\,2k+\tfrac{1}{2} < s < 2k+\tfrac{5}{2}.$$
We will assume that $p$ is odd, since under this assumption, \autoref{mHS} preserves its symmetry in the sense that $u_0(-x)=-u_0(x)$ implies that $u(t,-x)=-u(t,x)$ for any $t$ belonging to the interval of existence. Second, we consider an initial boundary value problem for the two-component Hunter-Saxton equation \eqref{2HS} on $(0,\frac{1}{2})$, with initial data belonging to $H^s(0,\frac{1}{2})\times H^{s-1}(0,\frac{1}{2})$ for $s\in[2,\tfrac{5}{2})$.

By a suitable transformation, we convert the initial boundary value problems on $(0,\frac{1}{2})$ into periodic Cauchy problems on the line and then apply results obtained in \cite{T05} and \cite{WW10,W10} to discuss local and global well-posedness and blow-up. Before we start, we provide some elementary results about odd periodic extensions of Sobolev functions defined on $(0,\frac{1}{2})$ in an introductory section and we recapitulate some results of \cite{T05,WW10,W10} for the modified periodic Hunter-Saxton equation and the periodic two-component Hunter-Saxton system.\\[-.25cm]

\textbf{Acknowledgement.} The financial support of the research training group GRK 1463 \emph{Analysis, Geometry and String Theory} (Leibniz University Hannover, Germany) is gratefully acknowledged. Cordial thanks also go to the referees whose comments helped to improve the initial version of the paper.
\section{Preliminaries}\label{sec_prelim}
In the theory of Sobolev spaces $W^{k}_m(\Omega)$ on a certain open domain $\Omega\varsubsetneq\R^n$ which does not behave too poorly, the following question is of interest: Is there an extension operator $E$ mapping functions defined on $\Omega$ to functions on $\R^n$ such that $Ex=x$ on $\Omega$ and such that $E\colon W^k_m(\Omega)\to W^{k}_m(\R^n)$ is continuous? For a wide range of results related to this problem, we refer the reader to \cite{A75}. In view of the purposes of this paper, we are in need of extensions in the case $n=1$ where $\Omega$ is the interval $(0,\frac{1}{2})$. We fix $m=2$ and ask for extensions of a Sobolev function defined on $(0,\frac{1}{2})$ to a function on the real axis; we also demand that the extended function is odd and periodic. Unfortunately, this case is not completely covered by \cite{A75} so that we will use extension results written down by Escher and Yin in \cite{EY09}, which will be very suitable for our purposes.
\defn We denote by $H^s(a,b)$ the $L_2$-Sobolev space of order $s\geq 0$ on the non-empty interval $(a,b)\subset\R$.
For $k\in\N_0$ and $2k+\frac{1}{2}<s< 2k+\frac{5}{2}$, let
$$D^s_k(a,b)=\set{v\in H^s(a,b)}{v^{(2k-2j)}(a)=v^{(2k-2j)}(b)=0,\;\forall j=0,\ldots,k}.$$
\enddefn
Observe that $D^s_0(a,b)$ is exactly the set of trace-zero functions in $H^s(a,b)$ for $s\in(\frac{1}{2},\frac{5}{2})$. Our approach to the modified Hunter-Saxton equation and the 2HS system is mainly based on an application of the following two lemmas.
\lem[\cite{EY09}; Lemma 2.12]\label{lem2.11} Let $s\in(\frac{1}{2},\frac{5}{2})$ be given and assume that $v\in H^s(0,\ell)$ with $v(0)=0$. Let furthermore
$$\tilde v(x)=
\left\{
\begin{array}{ll}
v(x),   & x\in[0,\ell), \\
-v(-x), & x\in(-\ell,0). \\
\end{array}
\right.
$$
Then $\tilde v\in H^s(-\ell,\ell)$.
\endlem\rm
\rem For $s\geq\frac{5}{2}$, under the assumption of Lemma~\ref{lem2.11}, one cannot conclude $\tilde v\in H^s(-\ell,\ell)$ generally. To understand this, it might be useful to recall that the norm in $H^{n+\sigma}(0,\ell)$, $n\in\N_0$, $\sigma\in(0,1)$, is given by
$$\Vert f\Vert_{H^{n+\sigma}(0,\ell)}^2=\Vert f\Vert^2_{H^n(0,\ell)}+\iint_{[0,\ell]^2}\frac{|f^{(n)}(x)-f^{(n)}(y)|^2}{|x-y|^{1+2\sigma}}\dx\d y$$
and that the estimate
$$\int_0^{\ell}\frac{|f(x)|^2}{x^{2s}}\dx\leq 2\left(1+\frac{4}{(2s-1)^2}\right)\iint_{[0,\ell]^2}\frac{|f(x)-f(y)|^2}{|x-y|^{1+2s}}\dx\d y,\quad f\in C_0^{\infty}(0,\ell),$$
for the double integral term is involved. The crucial point is that this estimate requires $s\in(\frac{1}{2},1)$. Combined with a density argument, this in turn only allows an odd periodic extension as in Lemma~\ref{lem2.11} for $s\in(\frac{1}{2},1)\cup[1,2]\cup(2,\frac{5}{2})$, cf.~\cite{EY09}.
\endrem
An extension result for Sobolev spaces of higher regularity is provided by the following lemma, in which we impose additional conditions on the function $v$.
\lem[\cite{EY09}; Lemma 2.13]\label{lem2.13} Let $v\in D^s_k(0,\ell)$ for $2k+\frac{1}{2}<s<2k+\frac{5}{2}$ with $k\in\N_0$ be given. Let furthermore
$$\tilde v(x)=
\left\{
\begin{array}{ll}
v(x),   & x\in[0,\ell), \\
-v(-x), & x\in(-\ell,0). \\
\end{array}
\right.
$$
Then $\tilde v\in D^s_k(-\ell,\ell)$.
\endlem\rm
For clarity, we recall the following: let $\S=\R/\Z$ denote the unit circle, i.e., $\S$ consists of equivalence classes of real numbers such that $x,y$ are equivalent if and only if $x-y\in\Z$. A representative system for this equivalence relation is the set $[0,1)$. A function $u\colon\S\to\R$ is defined on $\R$ and is periodic with period $1$.

We have the following local well-posedness result for the modified Hunter-Saxton equation on the circle.
\thm[\cite{T05,Y04}]\label{LWPTiglay}
Let $p\geq 1$ be any positive integer and $s>\frac{3}{2}$. Given the initial data $u_0\in H^s(\S)$, the periodic Cauchy problem for \autoref{mHS} has a unique solution belonging to
$$C\left([0,\delta); H^s(\S)\right)\cap C^1\left([0,\delta);H^{s-1}(\S)\right),$$
for some $\delta>0$, and this solution depends continuously on the initial data $u_0$. If $p=1$, there is a maximal interval of existence for the associated solution and the maximal existence time $T>0$ is independent of $s$ in the following sense: If
$$u=u(\cdot,u_0)\in C\left([0,T); H^s(\S)\right)\cap C^1\left([0,T);H^{s-1}(\S)\right)$$
and $u_0\in H^{s'}(\S)$ with $s'\neq s>\tfrac{3}{2}$, then
$$u\in C\left([0,T); H^{s'}(\S)\right)\cap C^1\left([0,T);H^{s'-1}(\S)\right),$$
with the same $T$. In particular, if $u_0\in H^{\infty}(\S)=\cap_{s\geq 0}H^s(\S)$, it follows that $u\in C([0,T);H^{\infty}(\S))$.
\endthm\rm
\rem Although it is an open problem whether \autoref{mHS} re-expresses geodesic motion on a diffeomorphism group of the circle for $p>1$, the above theorem is proved by rewriting \eqref{mHS} as an ordinary differential equation on a circle diffeomorphism group in \cite{T05}.
\endrem
Note that we allow $T=\infty$ in the above theorem. In this case, we call the corresponding solution of the Hunter-Saxton equation a global solution. If $T<\infty$, the corresponding solution is a local-in-time solution.

An analogue well-posedness result for the 2-component system \eqref{2HS} and the spaces $H^s(\S)\times H^{s-1}(\S)$ for $s\geq 2$ is proved in \cite{W10} where the author applies Kato's semigroup approach. The Hunter-Saxton system also allows for blow-up solutions; these are solutions $(u,\rho)$ with a finite existence time $T$ and the additional property that the $H^s\times H^{s-1}$-norm of $(u,\rho)$ is unbounded as $0<t\to T$ from below. A natural criterion for blow-up is that the first order spatial derivative of the finite-time solution $u$ is unbounded from below as $t$ approaches $T<\infty$ from below. This is motivated from corresponding results for the Camassa-Holm equation in \cite{CE00} where the blow-up is related to wave breaking in the sense that the slope of the wave profile tends to $-\infty$ near the wave crests as $t\to T$. We have the following precise blow-up result.
\thm[\cite{WW10}; Theorem 4.4]\label{thmblowup} Fix $(a,\kappa)\in\{2\}\times\R_+$. For any $(u_0,\rho_0)\in H^2(\S)\times H^1(\S)$ let $T$ denote the maximal existence time of the solution $(u,\rho)$ to the Hunter-Saxton system \eqref{2HS}. Then the solution $(u,\rho)$ blows up in finite time if and only if
$$\liminf_{t\to T^-}\inf_{x\in\S}u_x(t,x)=-\infty.$$
\endthm\rm
There are also some global well-posedness results for the Hunter-Saxton equation in the literature, cf., e.g., \cite{BC05,L07'',WW10}. For our purposes, we will consider the following persistence result for solutions of the 2HS system in $H^1(\S)\times H^0(\S)$.
\thm[\cite{W10}; Proposition 4.2]\label{thmglobal} Fix $(a,\kappa)\in\{2\}\times\R_+$. Let $(u_0,\rho_0)\in H^s(\S)\times H^{s-1}(\S)$ for $s\geq 2$ be given. Then the solution $(u,\rho)$ to the Hunter-Saxton system \eqref{2HS} will not blow-up in $H^1(\S)\times H^0(\S)$ in any finite time.
\endthm\rm
\section{The modified Hunter-Saxton equation on the interval $(0,\frac{1}{2})$}
Let us first study the following initial boundary value problem for the modified Hunter-Saxton equation on the finite interval $(0,\frac{1}{2})$:
\bea
\left\{
\begin{array}{rlll}
      \partial_tu+u^p\partial_xu & = & \frac{1}{2}\partial_x^{-1}\left(\partial_x(u^p)\partial_xu\right), & t>0,\;x\in(0,\frac{1}{2}), \\
      u(0,x) & = & u_0(x), & x\in[0,\frac{1}{2}], \\
      u(t,0) & = & 0, & t\geq 0, \\
      u(t,\frac{1}{2}) & = & 0, & t\geq 0.
\end{array}
\right.
\label{IVPmHS}
\eea
Counter to first intuition, the operator $\partial_x^{-1}$ is not the integral with upper boundary $x$; instead, we define
$$(\partial_x^{-1}u)(x)=\int_a^xu(y)\d y-\int_0^1\int_a^yu(z)\d z\d y,$$
for some $a\in[0,1)$, cf.~\cite{T05}. In our next lemma, we establish that for $p$ odd, a solution of \autoref{mHS} preserves its symmetry.
\lem\label{lem_sym} Let $u_0$ be an odd function and denote by $u$ the solution of \eqref{mHS} with $u(0,x)=u_0(x)$ for all $x\in\S$, obtained in Theorem~\ref{LWPTiglay}. Assume further that $p$ is odd. Then $u$ is an odd function.
\endlem\rm
\proof The initial value problem for \autoref{mHS} is equivalent to
\beq
\left\{
\begin{array}{rclll}
u_{tx}+\frac{1}{2}pu^{p-1}u_x^2+u^pu_{xx}&=&0, & \text{for} & t>0,\\
u &=& u_0, & \text{for} & t=0.
\end{array}
\right.
\label{mHS'}\eeq
Now since $p$ and $u_0$ are odd, it follows that the problem \eqref{mHS'} is invariant under the transformation $(u,x)\mapsto (-u,-x)$. Thus if $u(x)$ is a solution, then $v(x)=-u(-x)$ is a second solution. Unity of the solution implies that $u=v$ and this achieves the proof.
\endproof
Our first main results read as follows.
\thm\label{LWP1a} Assume that $u_0\in H^s(0,\frac{1}{2})$ with $\frac{3}{2}<s<\frac{5}{2}$ and $u_0(0)=u_0(\frac{1}{2})=0$ is given. Fix $p=1$. Then there exists a maximal $T=T(u_0)>0$ and a unique solution $u(t,x)$ to \autoref{IVPmHS} with initial value $u_0$ such that
$$u\in C\left([0,T);\,H^s(0,\tfrac{1}{2})\right)\cap C^1\left([0,T);\,H^{s-1}(0,\tfrac{1}{2})\right).$$
Moreover, the solution depends continuously on the initial data, i.e., the mapping $u_0\mapsto u(\cdot,u_0)$,
$$H^s(0,\tfrac{1}{2})\to C\left([0,T);\,H^s(0,\tfrac{1}{2})\right)\cap C^1\left([0,T);\,H^{s-1}(0,\tfrac{1}{2})\right)$$
is continuous, and the maximal time $T$ is independent of $s$: If
$$u(\cdot,u_0)\in C\left([0,T);\,H^s(0,\tfrac{1}{2})\right)\cap C^1\left([0,T);\,H^{s-1}(0,\tfrac{1}{2})\right)$$
is a solution and $u_0\in H^{s'}(\S)$ with $\frac{3}{2} < s' < \frac{5}{2}$, $s'\neq s$, then
$$u\in C\left([0,T);\,H^{s'}(0,\tfrac{1}{2})\right)\cap C^1\left([0,T);\,H^{s'-1}(0,\tfrac{1}{2})\right)$$
with the same $T$.
\endthm
\thm\label{LWP1b} Assume that $u_0\in H^s(0,\frac{1}{2})$ with $u_0(0)=u_0(\frac{1}{2})=0$ and $\frac{3}{2}<s<\frac{5}{2}$ is given and that $p$ is odd. Then there exists a positive number $\delta$ and a unique solution $u(t,x)$ to \autoref{IVPmHS} such that
$$u\in C\left([0,\delta);\,H^s(0,\tfrac{1}{2})\right)\cap C^1\left([0,\delta);\,H^{s-1}(0,\tfrac{1}{2})\right).$$
Moreover, the solution depends continuously on the initial data, i.e., the mapping $u_0\mapsto u(\cdot,u_0)$,
$$H^s(0,\tfrac{1}{2})\to C\left([0,\delta);\,H^s(0,\tfrac{1}{2})\right)\cap C^1\left([0,\delta);\,H^{s-1}(0,\tfrac{1}{2})\right)$$
is continuous.
\endthm\rm
We only give a proof of Theorem~\ref{LWP1a}. The key argument in the proof has already been used in \cite{EY09} where the authors study a similar initial boundary value problem for the rod equation \cite{D98,D00} and the $b$-equation.
\proof Our strategy is to convert the initial boundary value problem for the Hunter-Saxton equation on the interval $(0,\frac{1}{2})$ into the Cauchy problem for the periodic Hunter-Saxton equation with period $1$ on the line and to apply Theorem~\ref{LWPTiglay}. To this end, we extend the initial data $u_0(x)$ defined on the interval $[0,\tfrac{1}{2}]$ to an odd periodic function defined on the real axis:
$$\tilde u_0(x)=
\left\{
\begin{array}{lll}
  u_0(x),   & x\in[n,n+\tfrac{1}{2}], & n\in\Z, \\
  -u_0(-x), & x\in[n-\tfrac{1}{2},n], & n\in\Z.
\end{array}
\right.
$$
By Lemma~\ref{lem2.11}, the function $\tilde u_0$ is an element of $H^s(\S)$, $s\in(\frac{3}{2},\frac{5}{2})$. In fact, $\tilde u_0$ is continuous at the tie points and
\bea\lim_{x\to\tfrac{1}{2}-0}\frac{\dw\tilde u_0(x)}{\dw x} & = & \lim_{x\to\tfrac{1}{2}-0}\frac{\tilde u_0(x)-\tilde u_0(\tfrac{1}{2})}{x-\tfrac{1}{2}}\nonumber\\
&=& \lim_{x\to\tfrac{1}{2}-0}\frac{u_0(x)}{x-\tfrac{1}{2}}\nonumber\\
&=& \lim_{x\to-\tfrac{1}{2}+0}\frac{-u_0(-x)}{x+\tfrac{1}{2}}\nonumber\\
&=& \lim_{x\to-\tfrac{1}{2}+0}\frac{\tilde u_0(x)-\tilde u_0(-\tfrac{1}{2})}{x-(-\tfrac{1}{2})}\nonumber\\
&=& \lim_{x\to-\tfrac{1}{2}+0}\frac{\dw\tilde u_0(x)}{\dw x}\nonumber\\
&=& \lim_{x\to\tfrac{1}{2}+0}\frac{\dw\tilde u_0(x)}{\dw x}\nonumber
\eea
so that $\tilde u_0\in C^1(\S)$ (which is a natural consequence of Sobolev's embedding theorem).
We may thus convert the initial boundary value problem for the HS equation on the interval $(0,\frac{1}{2})$ into the following problem:
\bea
\left\{
\begin{array}{rlll}
      \partial_t\tilde u+\tilde u\partial_x\tilde u & = & \frac{1}{2}\partial_x^{-1}\left((\partial_x\tilde u)^2\right), & t>0,\;x\in\R, \\
      \tilde u(t,x)   & = & \tilde u(t,x+1), & t\geq 0,\;x\in\R, \\
      \tilde u(0,x)   & = & \tilde u_0(x), & x\in\R, \\
      \tilde u_0(-x) & = & -\tilde u_0(x), & x\in\R, \\
      \tilde u_0(0)   & = & 0, \\
      \tilde u_0(\tfrac{1}{2}) & = & 0. \\
\end{array}
\right.
\label{HStilde}
\eea
Theorem~\ref{LWPTiglay} ensures the existence of a unique solution $\tilde u$, defined on a maximal time interval $[0,T)$, such that
$$\tilde u\in C\left([0,T);\,H^s(\S)\right)\cap C^1\left([0,T);\,H^{s-1}(\S)\right)$$
and such that the map $\tilde u_0\mapsto \tilde u(\cdot,\tilde u_0)$,
\beq H^s(\S)\to C\left([0,T);\,H^s(\S)\right)\cap C^1\left([0,T);\,H^{s-1}(\S)\right)\label{regularity}\eeq
is continuous. Furthermore the maximal time $T$ does not depend on $s$. By Lemma~\ref{lem_sym} we know that the solution $\tilde u(t,\cdot)$ is odd for any $t\in[0,T)$ and hence $\tilde u(t,0)\equiv 0$ on $[0,T)$. A further application of Lemma~\ref{lem_sym} ensures that
$\tilde u(t,\tfrac{1}{2})=\tilde u(t,-\tfrac{1}{2})=-\tilde u(t,\tfrac{1}{2})$
and hence $\tilde u(t,\tfrac{1}{2})\equiv 0$ on $[0,T)$. Let $u$ denote the restriction of $\tilde u$ to $[0,T)\times[0,\frac{1}{2}]$. Then $u$ is a solution of the initial boundary value problem \eqref{IVPmHS} with $p=1$ and its maximal time of existence is independent of the Sobolev index. Since $\tilde u$ depends continuously on $\tilde u_0$, it is clear that $u$ depends continuously on $u_0$. To see that the solution $u$ to \eqref{IVPmHS} is unique, we assume that there is a second solution $v\neq u$ with initial value $u_0$. First we may extend $v$ to a function
$$
\tilde v(t,x)=
\left\{
\begin{array}{lll}
  v(t,x),   & x\in[n,n+\tfrac{1}{2}], & n\in\Z, \\
  -v(t,-x), & x\in[n-\tfrac{1}{2},n], & n\in\Z.
\end{array}
\right.
$$
Then $\tilde v$ is a solution of \eqref{HStilde} and unity of the solution implies that $\tilde v\equiv\tilde u$ and hence $v\equiv u$, which is a contradiction. This achieves the proof.
\endproof
Observe that our method of proof is not suitable to study more general solutions of class
$$C\left([0,T);\,H^s(0,\tfrac{1}{2})\right)\cap C^1\left([0,T);\,H^{s-1}(0,\tfrac{1}{2})\right),\quad s\geq\tfrac{5}{2},$$
but we may study the initial boundary value problem
\bea
\left\{
\begin{array}{rlll}
      \partial_tu+u^p\partial_xu & = & \frac{1}{2}\partial_x^{-1}\left(\partial_x(u^p)\partial_xu\right), & t>0,\;x\in(0,\frac{1}{2}), \\
      u(0,x) & = & u_0(x), & x\in[0,\frac{1}{2}], \\
      u^{(2k-2j)}(t,0) & = & 0, & t\geq 0,\, j=0,\ldots,k,\\
      u^{(2k-2j)}(t,\frac{1}{2}) & = & 0, & t\geq 0,\, j=0,\ldots,k,
\end{array}
\right.
\label{IVPmHS2}
\eea
for which we have, by means of Lemma~\ref{lem2.13}, the following local well-posedness result, cf.~also \cite{EY09}.
\thm\label{LWP2} Assume that $u_0\in D_k^s(0,\frac{1}{2})$, where $k\in\N$ and $2k+\frac{1}{2}<s<2k+\frac{5}{2}$, and that $p$ is odd. Then there exists a positive number $\delta$ and a unique solution $u(t,x)$ to \autoref{IVPmHS2} such that
$$u\in C\left([0,\delta);\,D_k^s(0,\tfrac{1}{2})\right)\cap C^1\left([0,\delta);\,D_k^{s-1}(0,\tfrac{1}{2})\right).$$
Moreover, the solution depends continuously on the initial data, i.e., the mapping
$$u_0\mapsto u,\quad D^s_k(0,\tfrac{1}{2})\to C\left([0,\delta);\,D_k^s(0,\tfrac{1}{2})\right)\cap C^1\left([0,\delta);\,D_k^{s-1}(0,\tfrac{1}{2})\right)$$
is continuous.
\endthm\rm
\proof The theorem can be proved in analogy to Theorem~\ref{LWP1a}. We first extend the inital data $u_0$ defined on $[0,\frac{1}{2}]$ to an odd and periodic function $\tilde u_0$ on the line. By Lemma~\ref{lem2.13}, $\tilde u_0\in D^s_k(-\tfrac{1}{2},\tfrac{1}{2})$. We then solve the associated periodic problem on the line whose solution is denoted as $\tilde u$. Since $p$ is odd, in view of Lemma~\ref{lem_sym}, $\tilde u$ and any of the derivatives $\tilde u^{(2k-2j)}$ is an odd function. The further conclusions follow as in the proof of Theorem~\ref{LWP1a}.
\endproof
In view of Theorem~\ref{LWPTiglay}, it is clear that for $p=1$ the above theorem holds true with a maximal existence time $T$.
\section{The two-component Hunter-Saxton system on the interval $(0,\tfrac{1}{2})$}
In this section, we concentrate our attention on the initial boundary value problem
\beq
\left\{
\begin{array}{rcll}
m_t                  & = & -m_xu-2u_xm-\kappa\rho\rho_x, & t>0,\;x\in(0,\tfrac{1}{2}),\\
m                    & = & -u_{xx},                      & t\geq 0,\;x\in[0,\tfrac{1}{2}],\\
\rho_t               & = & -(\rho u)_x,                  & t>0,\;x\in(0,\tfrac{1}{2}),\\
u                    & = & u_0,                          & t=0,\;x\in[0,\tfrac{1}{2}],\\
\rho                 & = & \rho_0,                       & t=0,\;x\in[0,\tfrac{1}{2}],\\
u(t,0)               & = & 0,                            & t\geq 0,\\
\rho(t,0)            & = & 0,                            & t\geq 0,\\
u(t,\tfrac{1}{2})    & = & 0,                            & t\geq 0,\\
\rho(t,\tfrac{1}{2}) & = & 0,                            & t\geq 0,
\end{array}
\right.
\label{IBVP2HS}\eeq
with $\kappa$ positive. Since \autoref{2HS} is invariant under the transformation $(u,\rho,x)\mapsto(-u,-\rho,-x)$, we have the following analog of Lemma~\ref{lem_sym}.
\lem Let $(u_0,\rho_0)$ be odd functions and denote by $(u,\rho)$ the solution of the two-component Hunter-Saxton equation with parameters $(a,\kappa)\in\{2\}\times\R_+$, $u(0,x)=u_0(x)$ and $\rho(0,x)=\rho_0(x)$ for all $x\in\S$. Then $u$ and $\rho$ are odd functions.
\endlem\rm
By the same arguments as in the proof of Theorem~\ref{LWP1a}, we find the following well-posedness result.
\thm\label{LWP2HS} Let $s\in [2,\tfrac{5}{2})$ and fix $(u_0,\rho_0)\in H^s(0,\tfrac{1}{2})\times H^{s-1}(0,\tfrac{1}{2})$ with $u_0(0)=\rho_0(0)=0$ and $u_0(\tfrac{1}{2})=\rho_0(\tfrac{1}{2})=0$. Then there is a maximal $T>0$ and a unique solution
$$(u,\rho)\in C\left([0,T);\;H^s(0,\tfrac{1}{2})\times H^{s-1}(0,\tfrac{1}{2})\right)\cap C^1\left([0,T);\;H^{s-1}(0,\tfrac{1}{2})\times H^{s-2}(0,\tfrac{1}{2})\right)$$
of the initial boundary value problem \eqref{IBVP2HS}, so that the mapping $(u_0,\rho_0)\mapsto (u,\rho)$
\bea
&H^{s}(0,\tfrac{1}{2})\times H^{s-1}(0,\tfrac{1}{2})\to C\left([0,T);\;H^s(0,\tfrac{1}{2})\times H^{s-1}(0,\tfrac{1}{2})\right)\nonumber\\
&\hspace{6cm}\cap C^1\left([0,T);\;H^{s-1}(0,\tfrac{1}{2})\times H^{s-2}(0,\tfrac{1}{2})\right)\nonumber
\eea
is continuous. Furthermore, the time $T$ does not depend on $s$ in the sense that if $s'\neq s\in[2,\tfrac{5}{2})$ and $(u_0,\rho_0)\in H^{s'}(0,\tfrac{1}{2})\times H^{s'-1}(0,\tfrac{1}{2})$, then
\bea&(u,\rho)\in C\left([0,T);\;H^{s'}(0,\tfrac{1}{2})\times H^{s'-1}(0,\tfrac{1}{2})\right)\nonumber\\
&\hspace{3cm}\cap C^1\left([0,T);\;H^{s'-1}(0,\tfrac{1}{2})\times H^{s'-2}(0,\tfrac{1}{2})\right)\nonumber\eea
with the same $T$.
\endthm\rm
Furthermore, we can derive the following blow-up result.
\thm Let $u_0\in H^2(0,\frac{1}{2})\times H^1(0,\tfrac{1}{2})$ with $u_0(0)=u_0(\frac{1}{2})=0$ and $\rho_0(0)=\rho_0(\frac{1}{2})=0$ be given. Then blow-up of the solution $(u,\rho)$ to \autoref{IBVP2HS} with $\kappa>0$ in finite time $T<\infty$ occurs if and only if
$$\liminf_{t\to T^-}\inf_{x\in[0,\frac{1}{2}]}u_x(t,x)=-\infty.$$
\endthm\rm
\proof As before, we extend the initial data $(u_0,\rho_0)$ defined on $[0,\frac{1}{2}]$ to odd and periodic functions $(\tilde u_0,\tilde\rho_0)$ on the line. Let $(\tilde u,\tilde\rho)$ denote the solution of the periodic 2HS equation with initial values $(\tilde u_0,\tilde\rho_0)$, obtained as in the proof of Theorem~\ref{LWP1a}. Recall that $(u,\rho)=(\tilde u,\tilde\rho)|_{[0,T)\times[0,\frac{1}{2}]}$ is the unique solution to \eqref{IBVP2HS} with initial values $(u_0,\rho_0)$. In view of Theorem~\ref{thmblowup} we conclude that $(\tilde u,\tilde\rho)$ blows up in finite time $T<\infty$ if and only if
$$\liminf_{t\to T^-}\inf_{x\in\S}\tilde u_x(t,x)=-\infty.$$
Since $\tilde u(t,\cdot)$ is odd, it follows that $\tilde u_x(t,\cdot)$ is even and hence
$$\liminf_{t\to T^-}\inf_{x\in\S}\tilde u_x(t,x)=\liminf_{t\to T^-}\inf_{x\in[0,\frac{1}{2}]}u_x(t,x)$$
so that our proof is completed.
\endproof
Global well-posedness of the Hunter-Saxton system on $(0,\frac{1}{2})$ holds in the following sense.
\thm Let $(u_0,\rho_0)\in H^s(0,\frac{1}{2})\times H^{s-1}(0,\tfrac{1}{2})$ for $s\in[2,\frac{5}{2})$ with $u_0(0)=u_0(\frac{1}{2})=0$ and $\rho_0(0)=\rho_0(\frac{1}{2})=0$ be given. Then the solution $(u,\rho)$ to \autoref{IBVP2HS} with $\kappa>0$ obtained in Theorem~\ref{LWP2HS} does not blow up in $H^1(0,\frac{1}{2})\times H^0(0,\frac{1}{2})$ for any finite time.
\endthm\rm
\proof Let $(\tilde u_0,\tilde\rho_0)$ be the odd and periodic extension of $(u_0,\rho_0)$. In view of Theorem~\ref{thmglobal}, we know that the associated odd and periodic solution $(\tilde u,\tilde\rho)$ of \eqref{2HS} stays bounded in $H^1(\S)\times H^0(\S)$. Moreover, $(u,\rho)=(\tilde u,\tilde\rho)$ confined to $[0,T)\times[0,\frac{1}{2}]$ is the unique solution to \eqref{IBVP2HS} with initial values $(u_0,\rho_0)$. This implies that
$$\Vert u\Vert_{H^1(0,\frac{1}{2})}^2=\int_0^{1/2}\left(u^2+u_x^2\right)\dw x\leq\int_0^{1}\left(\tilde u^2+\tilde u_x^2\right)\dw x=\Vert\tilde u\Vert_{H^1(\S)}^2<\infty$$
and
$$\Vert \rho\Vert_{H^0(0,\frac{1}{2})}^2=\int_0^{1/2}\rho^2\dx\leq\int_0^{1}\tilde\rho^2\dx=\Vert\tilde\rho\Vert_{H^0(\S)}^2<\infty$$
for any finite time.
\endproof
\rem Comparing the results of Theorem 5.1 and Theorem 5.3 in \cite{WW10}, we see that for $s\in(2,\frac{5}{2}]$ it seems to be difficult to obtain global well-posedness for the Hunter-Saxton system in $H^s\times H^{s-1}$ whereas there are results for $s=2$ (with an additional sign condition on $\rho_0$) and $s>5/2$. Recall that our local well-posedness result does not cover the case $s\geq 5/2$. That is why we only have a weaker persistence result for the solutions of the Hunter-Saxton system \eqref{IBVP2HS} in $H^1\times H^0$.
\endrem
\section{Open Problems}
It is an open problem to discuss the initial boundary value problem for the following variant of the Hunter-Saxton equation. Let us denote by $\mu$ the mean of a periodic function; in particular, for 1-periodic functions $u$, $\mu$ is the integral operator
$$\mu(u)=\int_0^1u(x)\dx.$$
Recall the discussion in the introductory section about director fields in nematic liquid crystals where we obtained the Hunter-Saxton equation from the variational principle. Assuming that there is a preferred direction of the director field (e.g., coming from an exterior magnetic field, acting on dipoles, which orients the fluid particles in the crystal) and retroactive action of other dipoles on a given director, the variational principle results in the follwing Euler-Lagrange equation
$$u_{tx}=-\frac{1}{2}u_x^2-uu_{xx}+2\left(u\mu(u)+\frac{1}{2}\mu(u^2)\right)$$
cf.~\cite{KLM08}. Differentiating this identity with respect to $x$ we find that
\beq -u_{txx}=-2\mu(u)u_x+2u_xu_{xx}+uu_{xxx}.\label{muHS}\eeq
Equation~\eqref{muHS} is called the $\mu$-Hunter-Saxton equation. Geometric aspects of this equation, the periodic Cauchy problem, traveling wave solutions and its integrable structure have been examined in \cite{KLM08,K10,LMT09}. As written down in \cite{KLM08}, \autoref{muHS} is well-posed in $H^s(\S)$ for $s>3/2$ and it is easy to see that \eqref{muHS} is invariant under the transformation $(u,x)\mapsto (-u,-x)$, which suggests to apply the arguments presented in the main body of the paper for an associated initial boundary value problem. However, there is the following difficulty: Any odd and periodic function has mean zero, so that the extension procedure applied in the proof of Theorem \ref{LWP1a} results in the HS-equation (and \emph{not} in its $\mu$-variant) for the extended function $\tilde u$, since $\mu(\tilde u)=0$.

Note however that there is a close relationship between the $\mu$HS equation and the Camassa-Holm equation; this can be seen conspicuously in the geometric picture: On the group of orientation-preserving diffeomorphisms of the circle, both equations, the CH and the $\mu$HS, re-express geodesic motion. Moreover, the geodesic flow is not only produced by an affine connection but arises also from a right-invariant metric on the diffeomorphism group which is compatible with the covariant derivative, \cite{CK02,KLM08}. Precisely, the inertia operator which generates the metric on the tangent space at the identity element is $1-\partial_x^2$ for the CH and $\mu-\partial_x^2$ for the $\mu$HS. The difference of the associated Green's functions is presented in Figure~3 of \cite{LMT09} and has values in the order of $10^{-3}$. Moreover, the results of \cite{CE00,KLM08,LMT09} show that both equations exhibit a similar behavior concerning well-posedness and blow-up. Thus, presumably, techniques used to study initial boundary value problems for the CH equation \cite{EY09} would also work for $\mu$HS.
\end{document}